\documentclass[11pt,leqno]{article}
\usepackage{amssymb,amsmath,amsthm}

\usepackage{enumerate}
\numberwithin{equation}{section}

\newtheorem{theorem}{Theorem}[section]
\newtheorem{proposition}[theorem]{Proposition}
\newtheorem{lemma}[theorem]{Lemma}

\theoremstyle{definition}
\newtheorem{definition}[theorem]{Definition}
\newtheorem{notation}[theorem]{Notation}
\newtheorem{example}[theorem]{Example}

\newtheorem{remark}[theorem]{Remark}

\newcommand{\C}{{\mathbb{C}}}
\newcommand{\N}{{\mathbb{N}}}

\newcommand{\Z}{{\mathbb{Z}}}
\newcommand{\cor}{{\bf{k}}}

\newcommand{\coro}{{\bf{k}}{\tiny (0)}}
\newcommand{\hcor}{{\bf{\widehat k}}}
\newcommand{\hcoro}{{\bf{\widehat k}}{\tiny (0)}}
\newcommand{\tens}[1][]{\mathbin{\otimes_{\raise1.5ex\hbox to-.1em{}#1}}}
\def\phi{{\varphi}}
\def\epsilon{\varepsilon}
\def\sha{\mathcal{A}}
\def\shf{\mathcal{F}}
\def\shi{\mathcal{I}}
\def\shl{\mathcal{L}}
\def\sho{\mathcal{O}}
\def\shp{\mathcal{P}}

\newcommand{\OO}{\sho}
\newcommand{\E}[1][]{\mathcal{E}_{#1}}
\newcommand{\HE}[1][]{\widehat{\mathcal{E}}_{#1}}
\newcommand{\W}[1][]{\mathcal{W}_{#1}}
\newcommand{\HW}[1][]{\widehat{\mathcal{W}}_{#1}}
\newcommand{\SW}[1][]{\mathcal{W}^s_{#1}}
\newcommand{\SHW}[1][]{\widehat{\mathcal{W}}^s_{#1}}
\newcommand{\TW}[1][]{\mathcal{W}^t_{#1}}
\newcommand{\TO}[1][]{\mathcal{O}^{\hbar}_{#1}}
\newcommand{\THO}[1][]{\widehat{\mathcal{O}}^{\hbar}_{#1}}
\newcommand{\TTO}[1][]{\mathcal{O}^{t,\hbar}_{#1}}
\newcommand{\STO}[1][]{\mathcal{O}^{s,\hbar}_{#1}}

\newcommand{\STHO}[1][]{\widehat{\mathcal{O}}^{s,\hbar}_{#1}}

\newcommand{\gr}{\mathop{\mathrm{gr}}}
\newcommand{\sts}{{\mathfrak{S}}}
\newcommand{\stx}{{\mathfrak{X}}}
\newcommand{\stm}{{\mathfrak{M}}}

\newcommand{\rmptt}{{\{\rm pt\}}}

\renewcommand{\to}[1][]{\xrightarrow[]{#1}}

\newcommand{\isoto}[1][]{\xrightarrow[#1]%
{{\raisebox{-.6ex}[0ex][-.6ex]{$\mspace{1mu}\sim\mspace{2mu}$}}}}

\newcommand{\sect}{\Gamma}
\newcommand{\rsect}{\mathrm{R}\Gamma}
\newcommand{\RD}{{\rm D}}
\newcommand{\Rb}{{\rm b}}
\newcommand{\oim}[1]{{#1}_*}
\newcommand{\eim}[1]{{#1}_!}
\newcommand{\reim}[1]{{R#1}_!}
\newcommand{\opb}[1]{#1^{-1}}
\newcommand{\tw}[1]{\widetilde{#1}}
\newcommand{\twh}[1]{\widehat{#1}}
\def\tot{{\rm tot}}
\newcommand{\res}{{\rm res}}
\newcommand{\md[1]}{{\rm Mod}({#1})}
\newcommand{\indlim}[1][]{\mathop{\varinjlim}\limits_{#1}}
\newcommand{\eqdot}{\mathbin{:=}}
\newcommand{\cl}{\colon}
\newcommand{\scbul}{\,\raise.4ex\hbox{$\scriptscriptstyle\bullet$}\,}
\newcommand{\FPI}{\shf\shp\shi}
\newcommand{\bnum}{\begin{enumerate}[{\rm(i)}]}
\newcommand{\enum}{\end{enumerate}}
\newcommand{\lp}{{\rm(}}
\newcommand{\rp}{{\rm)}}
\newcommand{\eq}{\begin{eqnarray}}
\newcommand{\eneq}{\end{eqnarray}}
\newcommand{\eqn}{\begin{eqnarray*}}
\newcommand{\eneqn}{\end{eqnarray*}}

\begin{document}

\title
{An algebra of deformation quantization for star-exponentials  on complex symplectic manifolds}
\date{July 9, 2006}
\author{Giuseppe Dito and Pierre Schapira}
\maketitle

\begin{abstract}
The cotangent bundle $T^*X$ to a complex manifold $X$ is 
classically endowed with
the sheaf of $\cor$-algebras $\W[T^*X]$ of deformation quantization,
where $\cor\eqdot \W[\rmptt]$ is a subfield of $\C[[\hbar,\opb{\hbar}]$. 
Here, we
construct a new sheaf of $\cor$-algebras $\TW[T^*X]$ which contains
$\W[T^*X]$ as a subalgebra and an extra central parameter $t$.  
We give the symbol calculus for this algebra and prove that  quantized
symplectic transformations operate on it. If   $P$  is any section of
order zero of $\W[T^*X]$, we show that 
$\exp(t\opb{\hbar} P)$ is well defined in $\TW[T^*X]$. 
\end{abstract}


\section*{Introduction}
A fundamental tool for spectral analysis in deformation quantization is the star-exponential \cite{BFFLS}. However,
at the formal level, the star-exponential does not make sense as a formal series in $\hbar$ and $\opb{\hbar}$. The goal
of this article is to construct a new sheaf of algebras on the cotangent bundle $T^*X$ to a complex manifold $X$ 
in which the star-exponential has a meaning
and such that quantized symplectic transformations operate on such algebras.

On the cotangent bundle $T^*X$ to a complex manifold $X$, there is a
well-known sheaf of filtered algebras called
deformation quantization  algebra by many authors (see \cite{BFFLS}, \cite{Ko}, etc.).
This algebra, denoted $\HW[T^*X]$ here,  is constructed in \cite{P-S}  as well
as its analytic counterpart $\W[T^*X]$.
The sheaf  $\W[T^*X]$ is 
similar to the sheaf $\E[T^*X]$ of microdifferential operators of
\cite{S-K-K}, but with an extra
central parameter ${\hbar}$, a substitute to the lack of homogeneity\footnote{
In this paper, we write $\E[T^*X]$ and $\W[T^*X]$  instead of the classical notations $\E[X]$ and $\W[X]$.}.
Here ${\hbar}$ belongs to the field $\cor\eqdot\W[\rmptt]$, a subfield of 
$\C[[\hbar,\opb{\hbar}]$. (Note that the notation $\tau=\hbar^{-1}$ is used in \cite{P-S}.)
When $X$ is affine and one denotes by
$(x;u)$ a point of  $T^*X$, a section $P$ of this sheaf 
on an open subset $U\subset T^*X$ is represented by its total symbol
$\sigma_\tot(P)=\sum_{-\infty< j\leq m}p_j(x;u)\hbar^{-j}$, with $m\in\Z$, 
$p_j\in\OO_{T^*X}(U)$, the $p_j$'s satisfying suitable inequalities and the
product being given by the Leibniz formula. 

In this paper, we construct a new sheaf of $\cor$-algebras
$\TW[T^*X]$, with an extra central holomorphic parameter $t$ defined in a
neighborhood of $t=0$, with the property that complex symplectic 
transformations may be
locally quantized as isomorphisms of algebras and there
are natural morphisms of $\cor$-algebras 
$\W[T^*X]\to[\iota]\TW[T^*X]\to[\res]\W[T^*X]$
whose composition is the identity on $\W[T^*X]$. 
We give the symbol
calculus on $\TW[T^*X]$, which extends naturally that of $\W[T^*X]$ 
(however, now we get series in $\hbar^{j}$ with $-\infty<j<\infty$)
and finally 
we show that, if $P$ is a section of $\W[T^*X]$ of order $0$,
then  $\exp(t\opb{\hbar} P)$ is well defined in $\TW[T^*X]$.
We also briefly discuss the case where $T^*X$ is replaced
with a general symplectic manifold.

Our construction is as follows. First, we add a central holomorphic parameter $s\in\C$
and consider the sheaf $\W[\C\times T^*X]$, the 
subsheaf of $\W[T^*(\C\times X)]$ consisting of sections not depending on 
$\partial_s$. Denoting by $a\cl \C\times T^*X\to T^*X$ the
projection, we first define an algebra 
$\SW[T^*X]\eqdot R^1\eim{a}\W[\C\times T^*X]$. The algebra
structure with respect to the $s$-variable  is given by 
convolution, as in the case
of the space $H^1_c(\C;\sho_\C)$. In order to replace this convolution 
product by an usual product, we define 
the sheaf $\TW[T^*X]$ as the ``formal'' Laplace transform with respect
to the variables $s\opb{\hbar}$ of the algebra $\SW[T^*X]$.

In a deformation quantization context, the existence of
$\exp(t\opb{\hbar} P)$  in $\TW[T^*X]$ gives a precise meaning to
the star-exponential \cite{BFFLS} of $P$ which is heuristically related
to the Feynman Path Integral of $P$.

\medskip\noindent
{\bf Acknowledgments.} We would like to thank 
Masaki Kashiwara for extremely useful conversations  and 
helpful insights. The first named author thanks Yoshiaki Maeda for warm hospitality
at Keio university where this work was finalized, and the JSPS for financial support.

\section{Symbols}\label{section:symbols}

\subsubsection*{The fields $\hcor$ and $\cor$}
We set $\hcor\eqdot \C[[\hbar,\opb{\hbar}]$.
Hence, an element $a\in \hcor$ is a series
\eqn
&&a= \sum_{-\infty<j\leq m}a_j\hbar^{-j}, \quad a_j\in\C,\quad m\in\Z.
\eneqn
Consider the following condition on $a$:
\eq\label{eq:defW0}
&&\left\{  \parbox{300 pt}{
there exist  positive 
constants $C,\epsilon$ such that 
$\vert a_{j}\vert \leq C\epsilon^{-j}(-j)!$ for all $j<0$.
}\right. 
\eneq
We denote by $\cor$ the subfield of $\hcor$ consisting of series 
satisfying \eqref{eq:defW0}.
\vskip2mm
\noindent{\bf Convention}~~We endow $\hcor$, hence  $\cor$, with
the  filtration associated to
\eq\label{ordh}
{\rm ord}(\hbar) = -1.
\eneq

The fields $\hcor$ and $\cor$ are 
$\Z$-filtered\footnote
{In the sequel, we shall say ``filtered'' instead of ``$\Z$-filtered''.} 
and contain the subrings $\hcoro$ and $\coro$, respectively. Note that 
 $\hcoro=\C[[\hbar]]$ and $\coro=\cor\cap\hcoro$.

\subsubsection*{The sheaves $\THO[X]$ and $\TO[X]$}
Let $(X,\sho_X)$ be a complex manifold.

\begin{definition}\label{def:thoi}
\bnum
\item
We denote by $\THO[X]$ the
sheaf $\sho_X[[\hbar,\opb{\hbar}]$. 
In other words, $\THO[X]$ is the filtered $\hcor$-algebra
defined as follows:
A section $f(x,{\hbar})$ of $\sho_X^{{\hbar}}$ of order $\leq m$ ($m\in\Z$) 
on an open set $U$ of $X$ is a series
\eq\label{eq:defshoxtau}
&&f(x,{\hbar})=\sum_{-\infty<j\leq m}f_j(x)\hbar^{-j},
\eneq
with $f_j\in \sho_X(U)$. 
\item
We denote by $\TO[X]$ the filtered $\cor$-subalgebra of $\THO[X]$ 
consisting of sections $f(x,{\hbar})$ as above satisfying:
\eq\label{eq:defW}
&&\left\{  \parbox{300 pt}{
for any compact subset $K$ of $U$ there exist  positive 
constants $C,\epsilon$ such that 
$\sup\limits_{K}\vert f_{j}\vert \leq C\epsilon^{-j}(-j)!$ for all $j<0$.
}\right. 
\eneq
\enum
\end{definition}
Note that
\eq\label{eq:shofilt}
\THO[X]\simeq \THO[X](0)\tens[\hcoro]\hcor,
&&\TO[X]\simeq \TO[X](0)\tens[\coro]\cor.
\eneq
(To be correct, we should have written $\cor_X$, the constant sheaf with values in $\cor$, 
instead of $\cor$ in these formulas, and similarly for $\coro$, $\hcoro$ and $\hcor$.)

Also note that there exist isomorphisms of sheaves (not of algebras)
\eq
&& \THO[X](0)\simeq \sho_{X\times\C}\hat{\vert}_{X\times\{0\}},
\label{eq:shohtau=shoh}\\
&& \TO[X](0)\simeq \sho_{X\times\C}\vert_{X\times\{0\}},\label{eq:shotau=sho}
\eneq
where $\sho_{X\times\C}\hat{\vert}_{X\times\{0\}}$ is the formal completion
of $\sho_{X\times\C}$ along the hypersurface $X\times \{0\}$ of $X\times\C$
and $\sho_{X\times\C}{\vert}_{X\times\{0\}}$ is the restriction 
of $\sho_{X\times\C}$ to  $X\times \{0\}$. 

Denoting by $t$ the coordinate on $\C$, the isomorphism \eqref{eq:shotau=sho} is given by the map
\eqn
&&\sho_X^{{\hbar}}(0)\ni\sum_{j\leq 0}f_j\hbar^{-j}
\mapsto \sum_{j\geq 0}f_{-j}\dfrac{t^{j}}{j!}\in
\sho_{X\times\C}\vert_{X\times\{0\}}.
\eneqn

\subsubsection*{The convolution algebra $H^1_c(\C;\sho_\C)$}
The results of  this  subsection are  well known and elementary. 
We recall them for  the reader's convenience.

We consider the complex line $\C$ endowed with a holomorphic coordinate $s$.
Using this coordinate, we identify the
sheaf $\sho_{\C}$ of holomorphic functions on $\C$ and the
sheaf $\Omega_{\C}$ of holomorphic forms on $\C$. 

The space $H^1_c(\C;\sho_{\C})$ is endowed with a 
structure of an algebra by
\eqn
H^1_c(\C;\sho_\C)
\times H^1_c(\C;\sho_\C)&\to& H^2_c(\C^2;\sho_{\C^2})\\
&\to& H^1_c(\C;\sho_\C),
\eneqn
where the first arrow is the cup product and the second arrow is the
integration along the fibers of the map $\C^2\to \C$,
$(s,s')\mapsto s+s'$. 

When representing the cohomology classes by holomorphic functions, the
convolution product is described as follows.

For a compact subset $K$ of $\C$, we identify the vector space 
$H^1_K(\C;\sho_\C)$ with the quotient space 
$\sect(\C\setminus K;\sho_\C)/\sect(\C;\sho_\C)$ and, if 
$f\in \sect(\C\setminus K;\sho_\C)$, we still denote by $f$ its image
in $H^1_K(\C;\sho_\C)$ or in $H^1_c(\C;\sho_\C)$. 
Let $K$ and $L$ be compact subsets of $\C$, let 
$f\in \sect(\C\setminus K;\sho_\C)$ and 
$g\in \sect(\C\setminus L;\sho_\C)$.  
 The convolution product $f* g$ is given by
\eq\label{eq:conv1}
&&f* g(z)=\frac{1}{2i\pi}\int_{\gamma}f(z-w)g(w)dw
\eneq
where $\gamma$ is a counter clockwise 
oriented circle which contains $L$ and $\vert z\vert$ is chosen big
enough so that $z+K$ is outside of the disc bounded by $\gamma$. 
It is an easy exercise to show  
that this definition does not depend on the
representatives $f$ and $g$, and that  
to interchange the role of $f$ and $g$ in the formula
\eqref{eq:conv1} modifies the 
result by a function defined all over $\C$, hence gives the same
result in $H^1_c(\C;\sho_{\C})$. Therefore, we obtain  a commutative 
algebra structure on $H^1_c(\C;\sho_{\C})$.

\begin{example}\label{exa:znstarzm}
\eqn
&&\frac{1}{z^{n+1}}* \frac{1}{z^{m+1}}
=\frac{(n+m)!}{n!m!}\frac{1}{z^{n+m+1}}.
\eneqn
\end{example}

\subsubsection*{The sheaf $\STO[X]$}

{}From now on, we shall concentrate our study on $\TO[X]$. 

\begin{notation}
We shall often denote by $\C_s$  the complex line $\C$
endowed with the coordinate $s$. 
\end{notation}

\begin{lemma}\label{le:Siu}
Let $Y$ be a complex  manifold and $Z$ a Stein submanifold of~$Y$. Then 
$H^j(Z;\TO[Y](0)\vert_Z)$ vanishes for $j\neq 0$.
\end{lemma}
\begin{proof}
Using the isomorphism \eqref{eq:shotau=sho}, 
we may replace the sheaf $\TO[Y](0)$ 
 with the sheaf $\sho_{Y\times\C_t}\vert_{t=0}$.
By a theorem of Siu \cite{Si}, $Z \times \{0\}$ 
admits a fundamental system of
open Stein neighborhoods in $Y\times\C_t$ and the result follows. 
\end{proof}

Let $X$ be a complex manifold. 
The manifold $\C_s\times X$ is thus endowed with the $\cor$-filtered sheaf 
$\TO[\C_s\times X]$.
Let $a\cl\C_s\times X\to X$ denote the projection.

\begin{lemma}\label{le:vanishcoh2}
\bnum
\item
One has the isomorphism
\eqn
&&R^j\eim{a}\TO[\C_s\times X]
\simeq R^j\eim{a}\TO[\C_s\times X](0)\tens[\coro]\cor.
\eneqn
\item
$R^j\eim{a}\TO[\C_s\times X](0)\simeq 0$ for $j\neq 1$.
\item
Let $U\subset \subset V\subset \subset W$ be three  open subsets
of $X$ and assume that $W$ is Stein. Then the natural morphism 
$\sect(W;R^1\eim{a}\TO[\C_s\times X])\to
\sect(U;R^1\eim{a}\TO[\C_s\times X])$
factorizes through 
\eqn
&&
\indlim[K\subset \C_s]
\sect((\C_s\setminus K)\times V;\TO[\C_s\times X])/
\sect(\C_s\times V;\TO[\C_s\times X]),
 \eneqn
where $K$ ranges over the family of compact subsets of $\C$.
\enum
\end{lemma}
\begin{proof}
(i) follows from the projection formula for sheaves 
({\em i.e.,} $\reim{a}(F\tens\opb{a}G)\simeq \reim{a}F\tens G$)
 and \eqref{eq:shofilt}.

\noindent
(ii) For $x\in X$, we have
\eqn
&&H^j(\reim{a}\TO[\C_s\times X](0))_x\simeq 
\indlim[K]H^j_K(\C_s\times \{x\};\TO[\C_s\times X](0)\vert_{\C_s\times\{x\}}).
\eneqn
Applying the distinguished triangle of functors
\eqn
&&\rsect_K(\C_s\times \{x\};\scbul)\to\rsect(\C_s\times \{x\};\scbul)\to 
\rsect((\C_s\setminus K)\times \{x\};\scbul)\to[+1]
\eneqn
to the sheaf $\TO[\C_s\times X](0)\vert_{\C_s\times\{x\}}$
we get the result by Lemma~\ref{le:Siu} for $j>1$ and the case $j=0$
follows from the principle of analytic continuation. 

\vspace{0.2cm}
\noindent
(iii) Recall first that if $W$ is a Stein manifold and if $W_1\subset\subset W$ is open, there exists a 
Stein open subset $W_2$ of $W$ with $W_1\subset\subset W_2\subset\subset W$.

For a compact subset $L$ of $X$,  
$\sect(L;R^1\eim{a}\TO[\C_s\times X])\simeq 
\sect(L;R^1\eim{a}\TO[\C_s\times X](0))\tens[\coro]\cor$. Hence, it is
enough to prove the result for $\TO[\C_s\times X](0)$.

By Lemma~\ref{le:Siu},
$H^j(D\times U;\TO[\C_s\times X](0))$ vanishes
for  $D$ open in $\C_s$, $U$ Stein open in $X$ and $j\neq 0$.
Therefore, 
$H^j_{K\times U}(\C_s\times U;\sho_{\C_s\times X}^{{\hbar}}(0))$ vanishes
for $j\neq 1$ and we get the exact sequence:
\eqn
&&0\to\sect(\C_s\times U;\TO[\C_s\times X](0))\to\sect((\C_s\setminus K)\times U;\TO[\C_s\times X](0))\\
&&\hspace{5cm}\to H^1_{K\times U}(\C_s\times U;\TO[\C_s\times X](0))\to 0.
\eneqn
\end{proof}

\begin{definition}\label{def:shostau}
We set $\STO[X]\eqdot R^1\eim{a}\TO[\C_s\times X]$.
\end{definition}

Clearly, $\STO[X]$ is a sheaf of filtered $\cor$-modules. 
By Lemma~\ref{le:vanishcoh2}, a section $f(s,x,{\hbar})$ of order $m$ 
of the sheaf $\STO[X]$ on a Stein open subset $W$ of $X$  
may be written on any relatively compact open subset
$U$ of $W$ as a series
\eqn
&&f(s,x,{\hbar})=\sum_{-\infty<j\leq m}f_j(s,x)\hbar^{-j},
\eneqn
where 
$f_j(s,x)$ is a holomorphic function on 
$(\C_s\setminus K_0)\times U$ for a compact set $K_0$ not depending on $j$ and
the $f_j$'s satisfy an estimate \eqref{eq:defW} on each compact subset $K$ of 
$(\C_s\setminus K_0)\times U$.

We shall extend the product \eqref{eq:conv1} to $\STO[X]$ as follows.
For two sections  
$f(s,x,{\hbar})=\sum_{-\infty<j\leq m}f_j(s,x)\hbar^{-j}$ and 
$g(s,x,{\hbar})=\sum_{-\infty<j\leq m'}g_j(s,x)\hbar^{-j}$ 
of $\STO[X]$, we set:
\eq\label{eq:conv2}
&&\left\{
\parbox{300 pt}{
$f(s,x,{\hbar})* g(s,x,{\hbar}) 
=\sum_{-\infty<j\leq m+m'}h_j(s,x)\hbar^{-j}$,\\
$h_k(s,x)=\sum_{i+j=k}
\frac{1}{2i\pi}\int_{\gamma}f_i(s-w,x)g_j(w,x)dw.$
}\right.
\eneq

\begin{proposition}
The sheaf $\STO[X]$ has a structure of a filtered commutative $\cor$-algebra.
\end{proposition}
\begin{proof} 
It  is easily checked that 
multiplication by $\opb{\hbar}$ induces an isomorphism of 
sheaves of $\cor$-modules $\STO[X](m)\isoto\STO[X](m+1)$. Hence  we just need to check
that the product of two sections of order $0$ is a section of order $0$.
Let $f(s,x,{\hbar})=\sum_{-\infty<i\leq 0}f_i(s,x)\hbar^{-i}$ and 
$g(s,x,{\hbar})=\sum_{-\infty<j\leq 0}g_j(s,x)\hbar^{-j}$ be in $\STO[X](0)$ and $K$
a compact subset of $(\C_s\setminus K_0)\times U$.
Let $\gamma$ be a counter clockwise 
oriented circle which contains $K_0$ and $s>R$  big
enough so that $s+K_0$ does not meet $\gamma$. Then for 
$w\in \gamma$ and $x\in K\cap (\C_s\setminus K_0)\times U$, we have:
$$
\vert \sum_{i+j=k,i,j\leq0} f_i(s-w,x)g_j(w,x)\vert\leq C^2 (-k)!\sum_{i+j=k,i,j\leq0} \epsilon^{-i-j}\frac{(-i)! (-j)!}{(-k)!}\leq 3 C^2  \epsilon^{-k} (-k)!.
$$
Hence $h(s,x,\hbar)=\sum_{-\infty<j\leq 0}h_k(s,x)\hbar^{-k}$ defined by \eqref{eq:conv2} is in $\STO[X](0)$.
\end{proof}

\subsubsection*{The Laplace transform and the algebra $\TTO[X]$}

In order to replace the convolution product in the $s$-variable with
the ordinary product, we shall apply a  kind of Laplace transform to
$\STO[X]$. 
\begin{definition}
On a complex manifold $X$, we denote by $\TTO[X]$ the
filtered sheaf of $\cor$-modules defined as follows. 
A section $f(t,x,{\hbar})$ of $\TTO[X](m)$ ({\em i.e.,} a section of order $m$) 
on an open set $U$ of $X$ is a series
\eq\label{eq:defshoxtau2}
&&f(t,x,{\hbar})
= \sum_{-\infty< j< \infty}f_j(t,x)\hbar^{-j}, 
\quad f_j\in \sect(U;\OO_{\C\times X\vert_{t=0}}),
\eneq 
with the condition that for any compact subset $K$ of $U$ there
exists $\eta>0$ such that  $f_{j}(t,x)$ is holomorphic
in a neighborhood of $\{\vert t\vert\leq \eta\}\times K$ and 
satisfies
\eq\label{eq:defTTO}
&&\left\{  \parbox{300 pt}{
there exist positive constants $C,\epsilon$ such that \\
$\sup\limits_{x\in K,\vert t\vert\leq\eta}\vert f_{j}(t,x)\vert 
\leq C\cdot\epsilon^{-j}(-j)!$ for all $j<0$,
}\right. 
\eneq
\eq\label{eq:defTTO+}
&&\left\{  \parbox{300 pt}{
 there exist positive 
constants $M$ and $R$ such 
that \\
$\sup\limits_{x\in K}\vert f_{j}(t,x)\vert
\leq M\dfrac{R^{j-m}}{(j-m)!}\vert t\vert^{j-m}$ for $\vert t\vert\leq\eta$
and all $j\geq m$.
}\right.
\eneq 
\end{definition}
Let $f(t,x,{\hbar})
= \sum_{-\infty< j< \infty}f_j(t,x)\hbar^{-j}$ and 
$g(t,x,{\hbar})
= \sum_{-\infty< j< \infty}g_j(t,x)\hbar^{-j}$ be two sections of 
$\TTO[X]$ of order $m$ and $m'$ respectively. Define formally
\eq\label{otproduct}
&&h(t,x,{\hbar})= \sum_{-\infty< j< \infty}h_j(t,x)\hbar^{-j},\quad
h_k(t,x)= \sum_{i+j=k}f_i(t,x)g_j(t,x).
\eneq
\begin{lemma}\label{le:TTOalgebra}
\bnum
\item 
Multiplication by $\opb{\hbar}$ induces an isomorphism of 
sheaves of $\coro$-modules $\TTO[X](m)\isoto\TTO[X](m+1)$.
\item
The product \eqref{otproduct} of a section  $f(t,x,{\hbar})\in\TTO[X](m)$
and a section $g(t,x,{\hbar})\in\TTO[X](m')$ is well defined and  belongs to 
$\TTO[X](m+m')$.
\enum
\end{lemma}
\begin{proof}
(i) 
\noindent (a) Let $f(t,x,{\hbar})= \sum_{-\infty< j< \infty}f_j(t,x)\hbar^{-j}\in\TTO[X](m)$, then
$\opb{\hbar}f(t,x,{\hbar})=  \sum_{-\infty< j< \infty}\tilde f_j(t,x)\hbar^{-j}$,
with $\tilde f_j = f_{j-1}$. For any integer $j<0$, we have:
\eqn
&&\sup\limits_{x\in K,\vert t\vert\leq\eta}\vert \tilde f_{j}(t,x)\vert 
= \sup\limits_{x\in K,\vert t\vert\leq\eta}\vert  f_{j-1}(t,x)\vert
\leq  C \epsilon^{-j+1} (-j+1)! \leq  (C\epsilon) (\epsilon e)^{-j} (-j)!.
\eneqn
Hence Condition~\eqref{eq:defTTO} is satisfied.

For $j\geq m+1$, we have:
\eqn
\sup\limits_{x\in K}\vert \tilde f_{j}(t,x)\vert 
= \sup\limits_{x\in K}\vert  f_{j-1}(t,x)\vert
&\leq & M \frac{R^{j-m-1}}{(j-m-1)!} \vert t\vert^{j-m-1},
\eneqn
which is simply Condition~\eqref{eq:defTTO+} for $m+1$ and 
$\opb{\hbar}f(t,x,{\hbar})\in \TTO[X](m+1)$.

\noindent (b) Let ${\hbar}f(t,x,{\hbar})=  \sum_{-\infty< j< \infty}\tilde f_j(t,x)\hbar^{-j}$,
with $\tilde f_j = f_{j+1}$. For any integer $j<-1$, we have:
\eqn
&&\sup\limits_{x\in K,\vert t\vert\leq\eta}\vert \tilde f_{j}(t,x)\vert 
= \sup\limits_{x\in K,\vert t\vert\leq\eta}\vert  f_{j+1}(t,x)\vert
\leq  C \epsilon^{-j-1} (-j-1)!
\leq  \frac{C}{\epsilon} \epsilon^{-j} (-j)!.
\eneqn
For $j=-1$, we have:
\eqn
&&\sup\limits_{x\in K,\vert t\vert\leq\eta}\vert \tilde f_{-1}(t,x)\vert 
= \sup\limits_{x\in K,\vert t\vert\leq\eta}\vert  f_{0}(t,x)\vert=A\geq 0, 
\eneqn
since $f_{0}(t,x)$ is holomorphic in a neighborhood of $\{\vert t\vert\leq \eta\}\times K$.
Set $C^\prime= \max\{\frac{A}{\epsilon},\frac{C}{\epsilon}\}$, then for all integer $j<0$, we have:
\eqn
\sup\limits_{x\in K,\vert t\vert\leq\eta}\vert \tilde f_{j}(t,x)\vert 
&\leq & C^\prime \epsilon^{-j} (-j)!,
\eneqn
and Condition~\eqref{eq:defTTO} is satisfied.

For $j\geq m-1$, we have:
\eqn
\sup\limits_{x\in K}\vert \tilde f_{j}(t,x)\vert 
= \sup\limits_{x\in K}\vert  f_{j+1}(t,x)\vert
&\leq & M \frac{R^{j-m+1}}{(j-m+1)!} \vert t\vert^{j-m+1},
\eneqn
which is Condition~\eqref{eq:defTTO+} for $m-1$ and ${\hbar}f(t,x,{\hbar})\in \TTO[X](m-1)$.
Therefore, multiplication by $\opb{\hbar}$ induces an isomorphism  
$\TTO[X](m)\isoto\TTO[X](m+1)$.

\vspace*{0.2cm}
\noindent (ii) By (i), we may assume $m=m'=0$. Let $f= \sum_{-\infty< i< \infty}f_i(t,x)\hbar^{-i}$
and $g= \sum_{-\infty< j< \infty}g_j(t,x)\hbar^{-j}$ be in $\TTO[X](0)$.
Let $K$ be a compact set. There exists $\eta>0$ such that
$f_i(t,x)$ and $g_j(t,x)$ are holomorphic  in a neighborhood of $\{\vert t\vert\leq \eta\}\times K$.
Conditions~\eqref{eq:defTTO} and \eqref{eq:defTTO+} guarantee the existence of the positive
constants $C_1$, $\epsilon_1$, $M_1$ and $R_1$ for the $f_i$'s, and 
$C_2$, $\epsilon_2$, $M_2$ and $R_2$ for the $g_j$'s. We set $C=\max\{C_1,C_2\}$, 
$\epsilon=\max\{\epsilon_1,\epsilon_2\}$, $M=\max\{M_1,M_2\}$ and $R=\max\{R_1,R_2\}$

We shall show that the product~\eqref{otproduct} is well defined.
Let $h_k(t,x)= \sum_{i+j=k}f_i(t,x)g_j(t,x)$. 

\noindent (a) Consider the case $k<0$. The sum defining $h_k$ can be divided into three parts:
\eq\label{hk}
h_k= \sum_{k<i<0}f_ig_{k-i} +  \sum_{i\geq0}f_i g_{k-i} + \sum_{j\geq 0}f_{k-j} g_j.
\eneq
The first sum is finite and defines a holomorphic function in a 
neighborhood of $\{\vert t\vert\leq \eta \}\times K$.

In the second sum, $k-i$ is strictly negative 
and for each term in this sum Conditions~\eqref{eq:defTTO} and \eqref{eq:defTTO+} give
the following estimates when  $x\in K$ and $\vert t\vert \leq \eta$:
\eqn
\vert f_i(t,x) g_{k-i}(t,x)\vert &\leq& M\frac{(R\eta)^i}{i!} C \epsilon^{i-k}(i-k)!\\
&\leq& C M \epsilon^{-k} (-k)! (R\eta\epsilon)^i \binom{i-k}{i}
\eneqn
Recall that $\sum_{i\geq0}\alpha^i \binom{n+i}{i} = \frac{1}{(1-\alpha)^{n+1}}$ for $\vert\alpha\vert<1$.
When $R\eta\epsilon<1$, $(R\eta\epsilon)^i \binom{i-k}{i}$ is the general term of an absolutely convergent
series. Let $\tilde\eta=\min\{\eta,\frac{1}{2R\epsilon}\}$. Then
the second sum in~\eqref{hk} converges uniformly on  $\{\vert t\vert\leq \tilde\eta \}\times K$.

The third sum is handled in a similar way and one gets the estimate:
\eqn
\vert f_{k-j}(t,x) g_{j}(t,x)\vert &\leq&  C  M  \epsilon^{-k} (-k)! (R\eta\epsilon)^j \binom{j-k}{j}.
\eneqn
It follows from the preceding that $h_k$ for $k<0$ is a holomorphic function in a 
neighborhood of $\{\vert t\vert\leq \tilde\eta \}\times K$.

Let now show that $h_k$ satisfies Condition~\eqref{eq:defTTO}.
For $x\in K$ and $\vert t\vert \leq \tilde\eta$, the first sum in~\eqref{hk} is bounded by:
\eqn
&&\vert \sum_{k<i<0}f_i(t,x)g_{k-i}(t,x)\vert \leq C^2\sum_{k<i<0} \epsilon^{-i}\epsilon^{i-k} (-i)!(i-k)!
\leq C^2 \epsilon^{-k}(-k)!.
\eneqn

For the second and third sums we have:
\eqn
\vert \sum_{i\geq0}  f_i(t,x) g_{k-i}(t,x)\vert 
&\leq&  CM  \epsilon^{-k} (-k)! \frac{1}{(1- R\tilde\eta\epsilon)^{-k+1}},\\
\vert \sum_{j\geq 0} f_{k-j}(t,x) g_{j}(t,x)\vert &\leq&  C  M 
\epsilon^{-k} (-k)! \frac{1}{(1-R\tilde\eta\epsilon)^{-k+1}}.
\eneqn
Let $\tilde \epsilon=\max\{\epsilon, \frac{\epsilon}{(1- R\tilde\eta\epsilon)}\}$.
For $x\in K$ and $\vert t\vert\leq \tilde\eta$, we find that:
\eqn
\vert h_k(t,x)\vert&\leq& (C^2 + \frac{2CM}{(1- R\tilde\eta\epsilon)})
\tilde\epsilon^{-k}(-k)!.
\eneqn
Hence $h_k$ satisfies Condition~\eqref{eq:defTTO}.

\vspace*{2mm}
\noindent (b) The case $k\geq0$. We again split
the sum defining $h_k$ into three parts:
\eq\label{hk+}
h_k= \sum_{0\leq i\leq k}f_ig_{k-i} +  \sum_{i<0}f_i g_{k-i} + \sum_{j< 0}f_{k-j} g_j.
\eneq
The first sum is a holomorphic function in a 
neighborhood of $\{\vert t\vert\leq \eta \}\times K$.

For each term in the second sum, we have the following estimates
when  $x\in K$ and $\vert t\vert \leq \eta$:
\eqn
\vert f_i(t,x) g_{k-i}(t,x)\vert \leq C \epsilon^{-i}(-i)! M\frac{R^{k-i}}{(k-i)!} \vert t\vert^{k-i}
\leq C M (\epsilon R \eta)^{-i} \frac{R^k}{k!} \vert t\vert^{k}.
\eneqn
$ (\epsilon R \eta)^{-i}$ is the general term of the geometric series, hence the second
sum in~\eqref{hk+} defines a holomorphic function in a 
neighborhood of $\{\vert t\vert\leq \tilde\eta \}\times K$ where 
$\tilde\eta=\min\{\eta,\frac{1}{2R\epsilon}\}$.

Similarly, for the third sum we have:
\eqn
\vert f_{k-j}(t,x) g_{j}(t,x)\vert 
&\leq& C M (\epsilon R \eta)^{-j} \frac{R^k}{k!} \vert t\vert^{k}.
\eneqn
Therefore $h_k$ for $k\geq0$ is a holomorphic function in a 
neighborhood of $\{\vert t\vert\leq \tilde\eta \}\times K$.

We now show that $h_k$ satisfies Condition~\eqref{eq:defTTO+} with $m=0$.
For $x\in K$ and $\vert t\vert \leq \tilde\eta$, the first sum in~\eqref{hk+} is bounded by:
\eqn
&&\vert \sum_{0\leq i \leq k}f_i(t,x)g_{k-i}(t,x)\vert 
\leq M^2 \frac{R^k}{k!} \vert t\vert^{k} \sum_{0\leq i \leq k}\binom{k}{i} 
\leq M^2 \frac{(2R)^k}{k!} \vert t\vert^{k}. 
\eneqn
For the second and third sums we find:
\eqn
\vert \sum_{i<0}  f_i(t,x) g_{k-i}(t,x)\vert 
&\leq&  C  M \frac{R \tilde\eta \epsilon}{1 - R\tilde\eta\epsilon} \frac{R^k}{k!} \vert t\vert^{k} \\
\vert \sum_{j< 0} f_{k-j}(t,x) g_{j}(t,x)\vert &\leq&  
CM \frac{R \tilde\eta \epsilon }{1 - R\tilde\eta\epsilon} \frac{R^k}{k!} \vert t\vert^{k} .
\eneqn
For $x\in K$ and $\vert t\vert\leq \tilde\eta$, we have:
\eqn
\vert h_k(t,x)\vert&\leq& (M^2 + 2 C  M \frac{R \tilde\eta \epsilon }{1 - R\tilde\eta\epsilon})
 \frac{(2R)^k}{k!} \vert t\vert^{k}.
\eneqn
Hence $h_k$ satisfies Condition~\eqref{eq:defTTO+} with $m=0$.

The product of $f\in\TTO[X](0)$ and $g\in\TTO[X](0)$ is well defined
and $fg \in\TTO[X](0)$.
\end{proof}

Therefore:

\begin{proposition}\label{pro:TTOalgebra}
The sheaf $\TTO[X]$ is naturally endowed with a structure of a
commutative filtered $\cor$-algebra.
\end{proposition}

Let $U$ be an open subset of $X$ and let 
$f(s,x,{\hbar})\in \sect((\C_s\setminus K)\times U;\TO[\C_s\times X])$. 
One defines formally the Laplace transform $\shl(f)$ of $f$ by
\eqn
&&\shl(f)(t,x,{\hbar})=
\frac{1}{2i\pi}\int_{\gamma}f(s,x,{\hbar})\exp(st\opb{\hbar})\ ds,
\eneqn
where $\gamma$ is a counter clockwise oriented circle centered at $0$
with radius $R\gg 0$.

\begin{example}
\eqn
\shl(s^{-n-1})= \hbar^{-n}t^n/n!,
&&\shl(\frac{1}{s-1})=\exp(t\opb{\hbar}).
\eneqn
\end{example}

\begin{lemma}
The Laplace transform induces a $\cor$-linear monomorphism
\eqn
&& \opb{s} \cdot \sho_X[[\opb{s}]][[\hbar,\opb{\hbar}]\hookrightarrow 
 \sho_X[[t]][[\hbar,\opb{\hbar}]]
\eneqn
\end{lemma}
\begin{proof}
One notices that the Laplace transform is given by:
\eqn
&&\sum_{-\infty<j\leq m}\sum_{n\geq 0}a_{n,j}s^{-n-1}\hbar^{-j}
\mapsto \sum_{j\leq m}\sum_{n\geq 0}\frac{a_{n,j}}{n!}t^{n}\hbar^{-n-j},
\eneqn
and the result follows.
\end{proof}

\begin{theorem}\label{th:laplacesymb1}
The Laplace transform induces a $\cor$-linear isomorphism of filtered
$\cor$-algebras
\eq\label{eq:laplacesymbol}
&&\shl\cl \STO[X]\isoto \TTO[X].
\eneq
\end{theorem}
\begin{proof}
(i) By Lemma~\ref{le:TTOalgebra}, it is enough to check that $\shl$
induces an isomorphism $\STO[X](0)\isoto \TTO[X](0)$.

\vspace{0.3cm}
\noindent
(ii) 
Let $W$ be a Stein open subset of $X$ and let $U$ be a 
relatively compact open subset of $W$. Let us develop a section 
$f(s,x,{\hbar})$ of  $\sect(W;R^1\eim{a}\sho_{\C_s\times X}^{{\hbar}}(0))$ 
with respect to $\opb{s}$ for $s>R$. We get
\eqn
\tilde f(s,x,{\hbar})&=&
\sum_{-\infty<j\leq 0}\tw f_j(s,x)\hbar^{-j}\\
&=&\sum_{-\infty<j\leq 0}\sum_{n\geq 0}f_{j,n}(x)s^{-n-1}\hbar^{-j}
\eneqn
with the following Cauchy's estimates:
\eqn
&&\left\{  \parbox{300 pt}{
for any compact subset $K$ of $U$
 there exist positive 
constants $C,\epsilon,R$ such that 
$\sup\limits_{x\in K}\vert f_{j,n}(x)\vert 
\leq C\epsilon^{-j}(-j)!R^{n}$.
}\right. 
\eneqn
Applying the Laplace transform to $\tilde  f(s,x,{\hbar})$ means to replace 
$s^{-n-1}$ with $\dfrac{t^n}{n!}\hbar^{-n}$. Hence, we find
\eqn
\shl(\tilde f)(t,x,{\hbar})&=& \sum_{-\infty<j<\infty}f_j(t,x)\hbar^{-j}=
\sum_{-\infty<j\leq 0}\sum_{n\geq 0}f_{j,n}(x)\frac{t^n}{n!}\hbar^{-j-n}
\eneqn
where
\eqn
f_j(t,x)&=&\sum_{j\leq n,0\leq n}f_{j-n,n}(x)\frac{t^n}{n!}
\eneqn
satisfies
\eqn
\vert f_j(t,x)\vert
&\leq&
C\sum_{j\leq n,0\leq n}\epsilon^{n-j}\frac{(n-j)!}{n!}
(\vert t\vert R)^n.
\eneqn
Let $\eta<\opb{(\epsilon R)}$.  It follows that  $f_j(t,x)$
is holomorphic in a neighborhood of $\{\vert t\vert\leq \eta\}\times K$.

Assume $j< 0$, $\vert t\vert\leq\eta$ and $x\in K$. We get
\eqn
&&\vert f_j(t,x)\vert
\leq C\epsilon^{-j}(-j)!\sum_{0\leq n}\frac{(n-j)!}{(-j)!n!}
(\eta\epsilon R)^n
\leq \frac{C}{1-\eta\epsilon R}(\frac{\epsilon}{1-\eta\epsilon R})^{-j}(-j)!,
\eneqn
hence Condition~\eqref{eq:defTTO} is satisfied.

Assume $j\geq 0$. We get for $\vert t\vert\leq\eta$ and $x\in K$ 
\eqn
&&\vert f_j(t,x)\vert \leq C\frac{\epsilon^{-j}}{j!}
\sum_{j\leq n}\frac{j!(n-j)!}{n!}
(\vert t\vert R\epsilon)^n
\leq \frac{C}{1-\eta\epsilon R}\frac{R^j}{j!}\vert t\vert ^{j},
\eneqn
hence Condition~\eqref{eq:defTTO+} for $m=0$ is satisfied 
and $\shl(f)(t,x,{\hbar})$ is in $\TTO[X](0)$.

\vspace{0.3cm}
\noindent
(iii) Conversely, let $f(t,x,{\hbar})$ be a section of $\TTO[X](0)$. 
We develop $f$ as
\eq\label{ftxh}
&&f(t,x,{\hbar})=\sum_{-\infty<j<\infty} f_j(t,x)\hbar^{-j}=\sum_{-\infty<j<\infty}\sum_{n\geq 0}n!f_{j,n}(x)
\frac{t^n}{n!}\hbar^{-n}\hbar^{-j+n}.
\eneq
For any compact set $K$, there exists $\eta>0$ such that
$ f_j(t,x)$ is holomorphic in a neighborhood of $\{\vert t\vert\leq \eta\}\times K$.
Conditions~\eqref{eq:defTTO} and ~\eqref{eq:defTTO+} give
the Cauchy's estimates
\eqn
&&\vert f_{j,n}(x)\vert\leq C\epsilon^{-j}(-j)!\eta^{-n}\text{ for }j< 0,\\
&&\vert f_{j,n}(x)\vert\leq M\frac{R^j}{j!}\eta^{j-n}\text{ for } j\geq 0.
\eneqn

Notice that Condition~\eqref{eq:defTTO+} for $j>0$ implies that
\eq\label{zero}
f_j(0,x)=\frac{\partial f_j}{\partial s}(0,x)=\cdots =
\frac{\partial^{j-1} f_j}{\partial s^{j-1}}(0,x)=0,
\eneq
or $f_{j,n}(x)=0$ for $0\leq n \leq j-1$.

The inverse Laplace transform consists formally in remplacing  $\dfrac{t^n}{n!}\hbar^{-n}$ by $s^{-n-1}$
in \eqref{ftxh}, we then 
get
$$
\opb{\shl}(f)(s,x,{\hbar})=\tilde f(s,x,{\hbar})=\sum_{-\infty<j< \infty}\sum_{n\geq 0}
n!f_{j+n,n}(x)s^{-n-1}\hbar^{-j}.
$$
Writing $\tilde f(s,x,{\hbar})=\sum_{-\infty<j< \infty} \tilde f_j(s,x) \hbar^{-j}$,
\eqref{zero} implies that $\tilde f_j(s,x)=0$ for $j\geq 1$.

Let $R_1>R$ large enough so that $(\eta R_1)^{-1}\leq 1$.
We shall check that the sum $\tilde f_j(s,x)= \sum_{n\geq 0} n!f_{j+n,n}(x)s^{-n-1}$ defines a holomorphic
function in a neighborhood of $\{\vert s\vert\geq R_1\}\times K$ for any $j\leq 0$.

For $j\leq 0$, let us split the sum $\tilde f_j(s,x)$ as
\eq\label{sumz}
\tilde f_j(s,x)= \sum_{n\geq -j} n!f_{j+n,n}(x)s^{-n-1} +  \sum_{0\leq n < -j} n!f_{j+n,n}(x)s^{-n-1}.
\eneq
In the first sum we have $n+j\leq 0$ and 
for $\vert s\vert \geq R_1$ and $x\in K$, we get from the Cauchy's estimates
\eq
\vert n!f_{j+n,n}(x)s^{-n-1}\vert &\leq&  n! M \frac{R^{n+j}}{(n+j)!}\eta^j \vert s\vert^{-n-1}
\leq \frac{M}{R_1} (\eta R)^{j}(-j)! \binom{n}{-j}(\frac{R}{R_1})^n.\label{esti} 
\eneq
The right-hand side is the general term of a convergent series since $R_1>R$ 
and we get the result by noticing that the second sum in \eqref{sumz} is finite.

Finally we shall show that $\tilde f_j(s,x)$ satisfies the required estimates.
From \eqref{esti}, the first sum in \eqref{sumz} is bounded by
\eqn
\vert\sum_{n\geq -j} n!f_{j+n,n}(x)s^{-n-1}\vert &\leq&
\frac{M}{R_1} (\eta R)^{j}(-j)!\sum_{n\geq -j} \binom{n}{-j}(\frac{R}{R_1})^n\\
&\leq& \frac{M}{R_1-R} \big(\frac{1}{\eta (R_1-R)}\big)^{-j}(-j)!.
\eneqn
Similarly, for the second sum we have
\eqn
\vert\sum_{0\leq n < -j} n!f_{j+n,n}(x)s^{-n-1}\vert &\leq&
\frac{C}{R_1}\epsilon^{-j}(-j)! \sum_{0\leq n < -j} \frac{1}{\binom{-j}{n}} (\frac{1}{\eta R_1})^n
\leq \frac{2C}{R_1}\epsilon^{-j}(-j)!,
\eneqn
where the last inequality follows from $(\eta R_1)^{-1}\leq 1$ and $ \sum_{0\leq n < -j} \frac{1}{\binom{-j}{n}}\leq 2$.

Combining these estimates we get for $j\leq0$
\eqn
\vert\tilde f_j(s,x)\vert \leq \tilde C \tilde \epsilon^{-j}(-j)!,
\eneqn
with $\tilde C=\max\{\frac{M}{R_1-R},\frac{2C}{R_1} \}$ and $\tilde \epsilon=\max\{\epsilon,\frac{1}{\eta (R_1-R)}\}$.

Therefore $\tilde  f(s,x,\hbar)=\sum_{j\leq0}\tilde f_j(s,x)\hbar^{-j}$ is a section of $\STO[X](0)$
and $\shl(\tilde f)(t,x,\hbar)= f(t,x,\hbar)$.

\vspace{0.3cm}
\noindent
(iv) The fact that $\shl$ is  a morphism of algebras follows easily
from Example~\ref{exa:znstarzm}. 
\end{proof}

\subsubsection*{The ring $\gr \TTO[X]$}

If $\sha$ is a filtered sheaf of rings, we denote as usual by $\gr \sha$ the
associated graded ring.

Let $\C_u$ be the complex line endowed with the coordinate $u$ and
denote by $b\cl X\times\C_u\to X$ the projection.

\begin{definition}\label{def:shoexp}
\bnum
\item
One denotes by $\sho_X^{\exp u}$ the subsheaf of $\C$-algebras
 on $X$ of the sheaf 
$\oim {b}\sho_{X\times\C_u}$ whose sections on an open set 
$U\subset X$ are the  holomorphic functions 
$f(x,u)$ on $U\times\C_u$ satisfying:
\eqn
&&\left\{  \parbox{300 pt}{
for any compact subset $K$ of $U$
 there exist positive 
constants $C,R$ such that 
$\sup\limits_{x\in K}\vert f(x,u)\vert \leq C\exp (R\vert u\vert)$.
}\right. 
\eneqn
\item
One sets  $\sho_X^{\exp t\opb{\hbar}}[\hbar,\opb{\hbar}]=\sho_X^{\exp t\opb{\hbar}}\tens[\C]\C[\hbar,\opb{\hbar}]$.
\enum
\end{definition}

\begin{proposition}
There is a natural isomorphism of graded sheaves of rings
\eqn
&& \gr\TTO[X] \simeq \sho_X^{\exp t \opb{\hbar}}[\hbar,\opb{\hbar}].
\eneqn
\end{proposition}
\begin{proof}
First note the isomorphism 
\eqn
&&\STO[X](0)/\STO[X](-1)\simeq  R^1\eim{a}\sho_{\C_s\times X},
\eneqn
from which we deduce the isomorphism
 \eqn
&&\gr\STO[X]\simeq  R^1\eim{a}\sho_{\C_s\times X}\tens[\C]\C[\hbar,\opb{\hbar}].
\eneqn  
The classical Paley-Wiener theorem says that the Laplace transform
induces an isomorphism between $H^1_c(\C;\sho_\C)$ and the space of entire
functions of exponential type.
An extension of this result with holomorphic parameters provides
an isomorphism
\eqn
&& \shl\cl R^1\eim{a}\sho_{\C_s\times X}\isoto
\sho_X^{\exp t \opb{\hbar} }
\eneqn
and the result follows.
\end{proof}

\subsubsection*{The formal case}

It is possible to replace $\TO[X]$ with $\THO[X]$ in 
the preceding constructions and to  set  
\eq
&&
\STHO[X]\eqdot R^1\eim{a}\THO[\C_s\times X].
\eneq
However the Laplace transform of $\STHO[X]$ does not seem to have
an easy description. Indeed, its sections are no longer germs of
holomorphic functions with respect to $t$ as shown in the next example.

\begin{example}
Consider a
sequence $\{c_{j}\}_{j\leq 0}$ of complex numbers and the section 
$f$ of $\STHO[X]$ given by 
\eqn
&&f(s,{\hbar})=\sum_{j\leq 0} \dfrac{c_{j}}{(s-1)}\hbar^{-j}.
\eneqn
Then, formally, the Laplace transform of $f$ is given by
\eqn
&&\shl(f)(t,{\hbar})=\sum_{j\leq 0}\sum_{n\geq 0}
c_{j}\frac{t^n}{n!}\hbar^{-n-j},
\eneqn and the coefficient of $\hbar^0$ is 
$\sum_{n\geq 0}c_{-n}\frac{t^n}{n!}$, which does not belong to 
$\sho_{\C_t}\vert_{t=0}$ in general. 
\end{example}

\section{The algebra $\W[T^*X]$}\label{section:W}

Let $(X,\sho_X)$ be a complex manifold. The cotangent bundle $T^*X$ is a 
homogeneous symplectic manifold
endowed with the $\C^\times$-conic sheaf of 
rings  $\E[T^*X]$ 
of finite-order microdifferential operators. This  ring is
filtered
and contains in particular the subring $\E[T^*X](0)$ of
operators of order $\leq 0$. 
This ring is constructed in
\cite{S-K-K} and we assume that the reader is familiar with this 
theory, referring to \cite{K3} or \cite{S} for an
exposition. 

On the symplectic manifold $T^*X$ there exists another (no more conic)
useful sheaf of rings constructed as follows (see \cite{P-S}). 
Let $\C$ be the complex line endowed with the coordinate $t$ and $(t;\tau)$ the associated
coordinates on $T^*\C$. 
Set $T^*_{\{\tau\neq0\}}(X\times\C) = \{(x,t;\xi,\tau);\tau\neq 0\}$
and consider the map
\eq\label{eq:projrho}
&&\rho\cl T^*_{\{\tau\neq0\}}(X\times\C)\to T^*X,\quad (x,t;\xi,\tau)\mapsto(x;\xi/\tau).
\eneq

Set
\eq\label{eq:Ehatt}
&&\E[T^*(X\times\C),\twh t]=\{P\in\E[T^*(X\times\C)];\ [P,\partial/\partial_t] = 0\}.
\eneq
The ring $\W[T^*X]$ 
 on $T^*X$ is given by
\eqn
&&\W[T^*X]\eqdot\oim\rho(\E[T^*(X\times\C),\twh t]).
\eneqn
In the sequel we set
\eq\label{hbartau}
&& \hbar\eqdot \opb{\tau} .
\eneq

The ring $\W[T^*X]$ is filtered and we denote by $\W[T^*X](j)$ the
subsheaf of $\W[T^*X]$ consisting of sections of order less or equal to $j$.  
The following result was obtained in~\cite{P-S}.

\begin{theorem}\label{th:W1}
\bnum
\item
The sheaf $\W[T^*X]$ is naturally endowed with a structure of a
filtered $\cor$-algebra and 
$\gr\W[T^*X]\simeq \sho_{T^*X}[\hbar,\opb{\hbar}]$.
\item 
Consider two complex manifolds $X$ and $Y$, two open subsets 
$U_X\subset  T^*X$ and $U_Y\subset T^*Y$ and a 
symplectic isomorphism $\psi:U_X\isoto U_Y$. Then, locally, $\psi$ may
be quantized as an isomorphism of filtered $\cor$-algebras 
$\Psi\cl\W[T^*X]\isoto\W[T^*Y]$ such that the isomorphism induced on the
graded algebras coincides with the isomorphism 
$\sho_{T^*X}[\hbar,\opb{\hbar}]\isoto 
\sho_{T^*Y}[\hbar,\opb{\hbar}]$ induced by $\psi$.
\enum
\end{theorem}

\subsubsection*{Total symbols}
Assume that $X$ is affine of dimension $n$, that is, 
 $X$ is open in some $\C$-vector space $V$ of dimension $n$.
\begin{theorem}
Assume $X$ is affine. There is an isomorphism of filtered sheaves of
$\cor$-modules \lp not of algebras\rp, called the ``total symbol'' morphism:
\eq\label{eq:totsymbmor}
&&\sigma_\tot\cl \W[T^*X]\isoto \sho^{{\hbar}}_{T^*X}.
\eneq
The total symbol of a product is given by the Leibniz formula.
Denote by $(x)$ a local coordinate system on $X$ and
denote by $(x,u)$ the associated  local symplectic coordinate system
on $T^*X$. 
If $Q$ is an operator of total symbol $\sigma_{\tot}(Q)$, then 
\eq\label{eq:leibniz1}
&&\sigma_{\tot}(P\circ Q)
=\sum_{\alpha\in\N^n} \dfrac{\hbar^{\vert\alpha\vert}}{\alpha !} 
\partial^{\alpha}_u\sigma_{\tot}(P)\cdot\partial^{\alpha}_x\sigma_{\tot}(Q).
\eneq 
\end{theorem}

The total symbol of a 
section $P\in \W[T^*X](U)$ is thus written as a formal series:
\eq\label{eq:totsymb}
&&\sigma_{\tot}(P) 
= \sum_{-\infty\leq j\leq m}p_j(x;u)\hbar^{-j}, 
\quad m\in\Z, \quad p_j\in\OO_{T^*X}(U),
\eneq
with the condition~\eqref{eq:defW}. 

Note that \eqref{eq:leibniz1} does not depend of the choice of a local
coordinate system on $X$  but only on the affine structure of $V$. 
Indeed, \eqref{eq:leibniz1} may be rewritten as 
\eqn
&&\sigma_{\tot}(P\circ Q)
= (\exp(\hbar\langle d_u,d_y\rangle)
\sigma_{\tot}(P)(x,u)\sigma_{\tot}(Q)(y,v))\vert_{x=y,u=v}.
\eneqn 
where 
$\langle d_u,d_y\rangle=\sum_{i=1}^n\partial_{u_i}\partial_{y_i}$ does
not depend on the affine coordinate system.

\begin{remark}
Let us identify $X$ with the zero section of $T^*X$. Then the sheaf
$\TO[X]$ (see Def.~\ref{def:thoi}) is isomorphic to the left coherent $\W[T^*X]$-module obtained as the
quotient of $\W[T^*X]$ by the left ideal generated the vector fields on 
$X$.
\end{remark}

\section{The algebra $\SW[T^*X]$}\label{section:SW}

\subsubsection*{Operations on $\W[]$}
Let $S$ be a complex manifold of complex dimension $d_S$. 
One defines 
the sheaf $\W[S\times T^*X]$ on $S\times T^*X$ as the
subsheaf of $\W[T^*(S\times X)]$ consisting of sections which commute with 
the holomorphic functions on $S$. Heuristically, 
$\W[S\times T^*X]$ is the sheaf $\W[T^*X]$ with holomorphic
parameters on $S$. 
For a morphism of complex manifolds $f\cl S\to Z$ 
we shall still denote by $f$ the
map $S\times X\to Z\times X$, as well as the map   
$S\times T^*X\to Z\times T^* X$. 
One denotes as usual by $\Omega_S$ the sheaf of holomorphic forms of
maximal degree and one sets for short:
\eq\label{eq:Wds}
&&\W[S\times T^* X]^{(d_S)}=\W[S\times
T^*X]\tens[\sho_S]\Omega_S.
\eneq

Let us recall well-known operations of the theory of
microdifferential operators. Although these results do not seem to be
explicitly written in the literature, their proofs are  straightforward
and will not be given here.

Let $f\cl S\to Z$ be a morphism of complex manifolds. 
The usual operations of inverse image 
$f^*\cl\opb{f}\sho_Z\to\sho_S$ 
and of direct image 
$\int_f\cl \reim{f}\Omega_S[d_S]\to \Omega_Z[d_Z]$
extend to $\W[S\times T^*X]$. More precisely, there exist 
morphisms of sheaves of $\cor$-modules (the second morphism holds in
the derived category $\RD^\Rb(\cor_{Z\times T^*X})$):
\eq
&&f^*\cl \opb{f}\W[Z\times T^*Z]
\to\W[S\times T^*X],\label{eq:invimW}\\
&&
\int_f\cl \reim{f}(\W[S\times T^*X]^{(d_S)}\,[d_S])\to 
\W[Z\times T^*X]^{(d_Z)}\,[d_Z],\label{eq:intSW}
\eneq
these morphisms having the following properties:
\begin{itemize}
\item 
they are functorial with respect to $f$, that is,
for a morphism of complex manifolds $g\cl Z\to W$, one has
$(g\circ f)^*\simeq f^*\circ g^*$ and
$\int_{g\circ f}=\int_g\circ\int_f$,  
and moreover the inverse (resp.\ direct) image of the
identity morphism is the identity,
\item
when $X$ is affine, $f^*$ and $\int_f$ commute with the total symbol
morphism \eqref{eq:totsymbmor}.
\end{itemize}

As a convention, we 
choose the morphism in \eqref{eq:intSW} so that the integral of
$\dfrac{ds}{s}\in H^1_c(\C_s;\Omega_{\C_s})$ is $1$. In other words,
\eqn
&&\int_a\frac{1}{s}=\frac{1}{2i\pi}\int_\gamma\frac{ds}{s},
\eneqn
where $\gamma$ is a counter clockwise oriented
circle around the origin.

\subsubsection*{The algebra $\SW[T^*X]$}
Denote by 
\eq\label{eq:mapa}
&&a\cl \C_s\times T^*X\to T^*X
\eneq
the projection. Then, after identifying the sheaves $\sho_{\C_s}$ and
$\Omega_{\C_s}$ by $f(s)\mapsto f(s)ds$, 
the sheaf $R^1\eim{a}\W[\C_s\times T^*X]$ is endowed with a
structure of a filtered $\cor$-algebra
by
\eqn
H^1_c(\C_s\times T^*X;\W[\C_s\times T^* X])
\times H^1_c(\C_{s'}\times T^*X;\W[\C_{s'}\times T^* X])
&&\\
&&\hspace{-95pt}
\to H^2_c(\C_{s,s'}^2\times T^*X;\W[\C^2_{s,s'}\times T^* X])\\
&&\hspace{-95pt}\to H^1_c(\C_s;\W[\C_s\times T^* X]),
\eneqn
where the first arrow is the cup product and the second arrow is the
integration along the fibers of the map $\C^2\to \C$,
$(s,s')\mapsto s+s'$. 

\begin{definition}\label{def:SW}
The sheaf $\SW[T^*X]$ of $\cor$-modules on $T^*X$ is given by
\eq\label{eq:defSW}
&&\SW[T^*X]= R^1\eim{a}(\W[\C_s\times T^* X]).
\eneq
\end{definition}
After identifying  the holomorphic function $\dfrac{1}{s}$ with the
cohomology class it defines in $H^1_c(\C_s;\sho_{\C_s})$,
we define the morphism of sheaves
\eq\label{eq:WtoSW}
&&\iota\cl  \W[T^*X]\to\SW[T^*X], \quad P\mapsto \dfrac{1}{s}P.
\eneq
Clearly, the morphism \eqref{eq:WtoSW} is a monomorphism of sheaves of 
$\cor$-algebras.

We define the morphism of sheaves
\eq\label{eq:SWtoW}
&&\res\cl \SW[T^*X]\to\W[T^*X]
\eneq
by  the integration morphism \eqref{eq:intSW} associated to 
the map \eqref{eq:mapa}.
Clearly, the morphism \eqref{eq:SWtoW} is a morphism of sheaves of 
$\cor$-algebras. Hence:

\begin{theorem}\label{th:convalg1}
\bnum
\item
The sheaf $\SW[T^*X]$ is naturally endowed with a structure of a
filtered $\cor$-algebra and 
$\gr\SW[T^*X]\simeq R^1\eim{a}\sho_{\C_s\times T^*X}[\hbar,\opb{\hbar}]$.
\item 
The monomorphism $\iota$ in \eqref{eq:WtoSW} is a 
morphism of filtered $\cor$-algebras,
the integration morphism $\res$ in \eqref{eq:SWtoW} is a morphism of
filtered $\cor$-algebras  and the composition 
$\res\circ\iota\cl \W[T^*X]\to \SW[T^*X]\to\W[T^*X]$  is the identity.
\item
Consider two complex manifolds $X$ and $Y$, two open subsets 
$U_X\subset  T^*X$ and $U_Y\subset T^*Y$ and a 
symplectic isomorphism $\psi:U_X\isoto U_Y$. Then, locally, $\psi$ may
be quantized as an isomorphism of filtered $\cor$-algebras 
$\Psi\cl\SW[T^*X]\isoto\SW[T^*Y]$ such that the isomorphism induced on the
graded algebras coincides with the isomorphism 
$R^1\eim{a}\sho_{\C_s\times T^*X}[\hbar,\opb{\hbar}]\isoto 
R^1\eim{a}\sho_{\C_s\times T^*Y}[\hbar,\opb{\hbar}]$ induced by $\psi$.
\item
Assume $X$ is affine. There is an isomorphism of filtered sheaves of
$\cor$-modules \lp not of algebras\rp, called the ``total symbol'' morphism:
\eq\label{eq:totsymbmorws}
&&\sigma_\tot\cl \SW[T^*X]\isoto \STO[T^*X].
\eneq
The total symbol of a product is given by the Leibniz formula with a convolution
product in the $s$ variable \lp see \eqref{eq:leibnizstar}\rp.
\enum
\end{theorem}
\begin{proof}
These results follow immediately from Theorem~\ref{th:W1}.
\end{proof}

Assume that $X$ is affine.
For each Stein open subset $W$ of $T^*X$ and each 
relatively compact open subset $U\subset\subset W$, 
a section $P$ of $\SW[T^*X]$ on  $W$ admits a 
total symbol
\eq\label{eq:totsymb3} 
\sigma_\tot(P)(s,x,u)&=& \sum_{-\infty< j\leq m}p_j(s,x;u)\hbar^{-j}, 
\quad m\in\Z
\eneq
where $p_j$ belongs to 
$\sect((\C_s\setminus K_0)\times   U;\sho_{\C_s\times T^*X})$,
for a compact subset $K_0$ of $\C_s$ 
which depends only on $P$ and $U$, and the $p_j$'s satisfy an estimate
as in \eqref{eq:defW} on each compact subset $K$ of 
$(\C_s\setminus K_0)\times U$.  

Consider now two sections $P$ and $Q$ of $\SW[T^*X]$ on a Stein open set $W$
with total symbols as in \eqref{eq:totsymb3} (replacing 
$p_j$ with $q_j$ and $m$ with $m'$ for $Q$).
Then the total symbol of $P\circ Q$ is given by the Leibniz formula:
\eq\label{eq:leibnizstar}
&&\sigma_{\tot}(P\circ Q)
=\sum_{\alpha\in\N^n} \dfrac{\hbar^{\vert\alpha\vert}}{\alpha !} 
\partial^{\alpha}_u\sigma_{\tot}(P)*\partial^{\alpha}_x\sigma_{\tot}(Q),
\eneq 
where, setting $f(s,x,u)=\partial^{\alpha}_u\sigma_{\tot}(P)(s,x;u)$
and $g(s,x,u)=\partial^{\alpha}_x\sigma_{\tot}(Q)(s,x;u)$, 
the product $f* g$ is given by \eqref{eq:conv2}.

\section{The Laplace transform and the algebra $\TW[T^*X]$}\label{section:TW}

The filtered $\cor$-algebra $\TW[T^*X]$ on $T^*X$ is the
algebra $\SW[T^*X]$, but with a different symbol calculus.

\begin{definition}\label{def:TW}
We set $\TW[T^*X]\eqdot \SW[T^*X]$. For $X$ affine,
 the total symbol morphism of $\cor$-modules (not of
algebras) 
\eq\label{eq:totsymb4} 
&&\sigma_\tot\cl \TW[T^*X]\isoto \TTO[T^*X]
\eneq
is the composition 
$\SW[T^*X]\isoto[\sigma_\tot]\STO[T^*X] \isoto[\shl]\TTO[T^*X]$.
\end{definition}

For $P$ a section of $\TW[T^*X]$ on a Stein open subset $V$ of
$T^*X$ and  an open subset $U\subset\subset V$, $\sigma_{\tot}(P)$ is written as a series
\eqn
&&\sigma_{\tot}(P)(t,x,u,{\hbar})
= \sum_{-\infty< j< \infty}p_j(t,x,u)\hbar^{-j}, 
\quad p_j\in\OO_{\C\times T^*X\vert_{t=0}}(U)
\eneqn
satisfying~\eqref{eq:defTTO} and~\eqref{eq:defTTO+}.

Applying Theorem~\ref{th:convalg1}, we get:

\begin{theorem}\label{th:1}
\bnum
\item
$\TW[T^*X]$ is a filtered $\cor$-algebra and 
$\gr\TW[T^*X]\simeq \sho_{T^*X}^{\exp t\opb{\hbar}}[\hbar,\opb{\hbar}]$
\lp see Definition~\ref{def:shoexp}\rp.
\item
The  morphism $\iota$ in \eqref{eq:WtoSW} induces a monomorphism of
filtered $\cor$-algebras $\iota\cl\W[T^*X]\hookrightarrow \TW[T^*X]$, 
the morphism $\res$ in \eqref{eq:SWtoW} 
induces a morphism of filtered $\cor$-algebras  $\res\cl \TW[T^*X]\to \W[T^*X]$
and the composition 
$\W[T^*X]\to \TW[T^*X]\to\W[T^*X]$  is the identity.
\item
Consider two complex manifolds $X$ and $Y$, two open subsets 
$U_X\subset  T^*X$ and $U_Y\subset T^*Y$ and a 
symplectic isomorphism $\psi:U_X\isoto U_Y$. Then, locally, $\psi$ may
be quantized as an isomorphism of filtered $\cor$-algebras 
$\Psi\cl\TW[T^*X]\isoto\TW[T^*Y]$ such that the isomorphism induced on the
graded algebras coincides with the isomorphism 
$\sho_{T^*X}^{\exp t\opb{\hbar}}[\hbar,\opb{\hbar}]\isoto 
\sho_{T^*Y}^{\exp t\opb{\hbar}}[\hbar,\opb{\hbar}]$ induced by $\psi$.
\item
Assume $X$ is affine. There is an isomorphism of filtered sheaves of
$\cor$-modules \lp not of algebras\rp, called the ``total symbol'' morphism:
\eq\label{eq:totsymbmorws}
&&\sigma_\tot\cl \TW[T^*X]\isoto \TTO[T^*X].
\eneq
The total symbol of a product is given by the Leibniz formula. 
\enum
\end{theorem}

For $P$ and $Q$ two sections of $\TW[T^*X]$ on an open subset $U$ of
$T^*X$, with $X$ affine, the total symbol of $P\circ Q$ is thus given by the 
formula:
\eq\label{eq:leibniz3}
&&\sigma_{\tot}(P\circ Q)
=\sum_{\alpha\in\N^n} \dfrac{\hbar^{\vert\alpha\vert}}{\alpha !} 
\partial^{\alpha}_u\sigma_{\tot}(P)\cdot\partial^{\alpha}_x\sigma_{\tot}(Q),
\eneq 
where the product 
$\partial^{\alpha}_u\sigma_{\tot}(P)\cdot\partial^{\alpha}_x\sigma_{\tot}(Q)$ 
is given by the usual commutative algebra structure of 
$\TTO[T^*X]$ of Lemma~\ref{le:TTOalgebra}.

\begin{remark}
In  Theorem~\ref{th:1}, the monomorphism $\W[T^*X]\to\TW[T^*X]$
is given on symbols by $\sigma_\tot(P)\mapsto \sigma_\tot(P)$ 
 and the morphism $\TW[T^*X]\to\W[T^*X]$
is given on symbols by 
$\sigma_\tot(P)(t,x;u,{\hbar})\mapsto \sigma_\tot(P)(0,x;u,{\hbar})$. 
\end{remark}

\subsubsection*{The formal case}
The above constructions also work when replacing the 
sheaf $\W[T^*X]$ with its formal counterpart, the sheaf $\HW[T^*X]$. Let us
briefly explain it.

Let $X$ be a complex manifold, as above. 
Replacing the sheaf of rings $\E[T^*X]$ on $T^*X$ with 
the sheaf of rings $\HE[T^*X]$ 
of formal microdifferential operators and proceeding as for $\W[T^*X]$,
we get the sheaf of rings $\HW[T^*X]$ 
of finite-order formal WKB-operators on $T^*X$. It is defined by
\eqn
&&\HW[T^*X]\eqdot\oim\rho(\HE[T^*(X\times\C),\twh t]).
\eneqn
When $X$ is affine  of dimension $n$, the total symbol morphism
induces an isomorphism of $\hcor$-modules
\eqn
&&\sigma_\tot\cl \HW[T^*X]\isoto \THO[T^*X]
\eneqn
and the symbol $\sigma_\tot(P\circ Q)$ is given by the Leibniz 
formula~\eqref{eq:leibniz1}.
Then by a similar construction as for $\SW[T^*X]$ we
construct the filtered sheaf of $\hcor$-algebras $\SHW[T^*X]$. Namely, we set
\eqn
&&\SHW[T^*X]\eqdot R^1\eim{a}\HW[\C\times T^*X].
\eneqn
If $X$ is affine, the total symbol morphism induces an isomorphism of
$\hcor$-modules $\SHW[T^*X]\isoto \STHO[T^*X]$ and the product is again
given by the Leibniz formula \eqref{eq:leibnizstar}.

However, as already noticed, the Laplace transform does not seem to behave as well for 
the formal case as for the analytic case, and we shall not construct the Laplace transform of 
$\STHO[T^*X]$.

\section{Remark: The algebra $\SW[\stx]$ on a symplectic manifold $\stx$}

\subsubsection*{The complex case}

Consider a complex symplectic manifold $\stx$. There exists an open covering 
$\stx=\bigcup_iU_i$ and complex symplectic isomorphisms $\phi_i\cl U_i\isoto V_i$
where the $V_i$'s are open in some cotangent bundles $T^*X_i$ of complex manifolds $X_i$. Set 
$\W[U_i]\eqdot\opb{\phi_i}\W[T^*X_i\vert_{V_i}]$. In general, the $\W[U_i]$'s do not glue 
in order to give a globally defined sheaf of algebras $\W[\stx]$ on $\stx$.
However the prestack $\sts$ on $\stx$ (roughly speaking, a prestack is a sheaf of categories) 
$U_i\mapsto \md[{\W[U_i]}]$ is a stack and the category $\md[{\W[\stx]}]\eqdot \sts(\stx)$
is well defined. Moreover, one can give a precise meaning to $\W[\stx]$ by replacing the notion of a
sheaf of algebras with that of an algebroid. We refer to \cite{K2} for the 
construction of (an analogue of) this stack in the contact complex case 
and to \cite{Ko} in the symplectic complex case for $\HW[\stx]$ and for the definition of an 
algebroid.
See also \cite{P-S} for a construction of $\W[\stx]$ (by a different method).
By adapting the construction of \cite{P-S}, one easily constructs the algebroid $\SW[\stx]$ associated 
with the locally defined sheaves of algebras $\SW[U_i]$. Details are left to the reader.

\subsubsection*{The real case}

Let $M$ be a real analytic manifold, $X$ a complexification of $M$
and denote by $\omega_X$ the canonical $2$-form on $T^*X$. The conormal bundle $T^*_MX$ 
is Lagrangian for ${\rm Re}\ \omega_X$ and symplectic for ${\rm Im}\ \omega_X$. In particular, the real manifold 
$T^*_MX$ is symplectic. For an open subset $U$ of $T^*_MX$, we set 
$\W[U]\eqdot\W[T^*X]\vert_{U}$.

Now, consider a real analytic symplectic manifold $\stm$. It is well known that 
it is possible to construct a globally defined sheaf of algebras $\W[\stm]$ on $\stm$
such that:
\begin{itemize}
\item
 there exists an open covering 
$\stm=\bigcup_{i\in I}U_i$ and real symplectic isomorphisms $\phi_i\cl U_i\isoto V_i$
where the $V_i$'s are open in the conormal bundles $T_{M_i}^*X_i$ for some real manifolds $M_i$
with complexification $X_i$,
\item 
$\W[\stm]\vert_{V_i}\simeq\opb{\phi_i}\W[{V_i}]$ for all $i\in I$. 
\end{itemize}
Replacing  $\stm$ with  $\C_s\times\stm$, one easily contructs the sheaf of algebras 
$\W[\C_s\times\stm]$ of sections with holomorphic parameter $s\in\C_s$. Setting
\eqn
&&\SW[\stm]\eqdot R^1\eim{a}\W[\C_s\times\stm]
\eneqn
we get a filtered $\cor$-algebra similar to the algebra $\SW[T^*X]$ of Definition~\ref{def:SW}.
Then, if $P$ belongs to $\W[\stm]$ and has order $0$, the section 
$\dfrac{1}{s-P}$ is well defined in $\SW[\stm]$.

\section{Applications}

As an application, let us construct the exponential of 
sections of order $0$ of $\W[T^*X]$. 

Consider a section $P$ of $\W[T^*X](0)$ on an open subset $U$ of $T^*X$. 
For each compact subset $K$
of $U$, there exists $R>0$ such that the section $s-P$ of 
$\SW[T^*X]$ defined on $\C_s\times U$ is invertible on 
$(\C_s\setminus D(0,R))\times K$, where   $D(0,R)$ denotes the closed disc
centered at $0$ with radius $R$.  Therefore
$\dfrac{1}{s-P}$ defines an element of 
$H^1_c(\C_s\times U;\W[\C_s\times T^* X])$, 
hence, an element of $\sect(U;\SW[T^*X])$. We still denote this section of 
$\SW[T^*X]$ on $U$ by $\dfrac{1}{s-P}$.

By developing $\dfrac{1}{s-P}$ as 
$\sum_{n\geq 0}\dfrac{P^n}{s^{n+1}}$
and applying the Laplace transform, we get formally:
$\shl(\frac{1}{s-P})=\exp(t\opb{\hbar} P)$.

\begin{notation}
We denote by $\exp(t\opb{\hbar} P)$ the image in $\TW[T^*X]$ of
the section $\dfrac{1}{s-P}$ of $\SW[T^*X]$.
\end{notation}

\begin{proposition}\label{expstar}
For $P\in\W[T^*X](0)$, there is a 
 section $\exp(t\opb{\hbar} P)\in \TW[T^*X]$ 
such that, when $X$ is affine: 
\eqn
&&\sigma_\tot(\exp(t\opb{\hbar} P))
       =\sum_{n\geq 0}\dfrac{(t\opb{\hbar} \sigma_\tot(P))^{\star n}}{n!},
\eneqn
where the star-product $f^{\star n}$ means the product given by the 
Leibniz formula~\eqref{eq:leibniz1}.
\end{proposition}

\begin{remark}
The Leibniz formula~\eqref{eq:leibniz1} is nothing but the standard or
normal or Wick star-product and Proposition~\ref{expstar} tells us that
the star-exponential \cite{BFFLS} of $P$ makes sense in $\TW[T^*X]$.
\end{remark}

In a holomorphic deformation quantization context, 
the star-exponential of $P$ is heuristically related
to the Feynman Path Integral  $\FPI(P)$ of $P$. Indeed, the Feynman
Path Integral of a Hamiltonian $H$ is the  
symbol of the evolution operator associated to $H$, the precise
relation being given  (see \cite{JD}) by 
\eqn
&&\exp(- xu\opb{\hbar}) \FPI(P) = \sigma_{\tot}(\exp(t\opb{\hbar} P)).
\eneqn

\begin{example}
As a simple example, take $X=\C$ and $P\in\W[T^*X](0)$ with 
 $\sigma_{\tot}(P) = p_0(t,x;u)= \theta x u$, $\theta\in \C$.
Up to a change of holomorphic symplectic coordinates, 
$\sigma_{\tot}(P)$ represents
the Hamiltonian of the harmonic oscillator in the holomorphic representation. 
Clearly $P$ is in $\W[\C](0)$, and the total symbol of $\exp(t\opb{\hbar} P)$ 
is easily computed:

\begin{align*}  
\frac{\partial}{\partial t}\sigma_{\tot}(\exp(t\opb{\hbar} P))
=&\sigma_{\tot}(\opb{\hbar} P\circ \exp(t\opb{\hbar} P))\\
 =&\opb{\hbar} \big(\sigma_{\tot}(P)\sigma_{\tot}(\exp(t\opb{\hbar} P)) + 
 {\hbar}\frac{\partial}{\partial u}\sigma_{\tot}(P) 
\frac{\partial}{\partial x}\sigma_{\tot}(\exp(t\opb{\hbar} P))\big)\\
 =& \opb{\hbar}\theta ux \sigma_{\tot}(\exp(t\opb{\hbar} P)) + 
\theta x \frac{\partial}{\partial x}\sigma_{\tot}(\exp(t\opb{\hbar} P)).
\end{align*}
Since $\sigma_{\tot}(\exp(t\opb{\hbar} P))|_{t=0} = 1$, the solution to 
the preceding equation is:
$$
\sigma_{\tot}(\exp(t\opb{\hbar} P)) = \exp\big((\exp(\theta t) -1) xu\opb{\hbar} \big). 
$$

The Feynman Path Integral
for the harmonic oscillator is well known in the Physics literature 
and is given by $\exp\big(\exp(\theta t) xu\opb{\hbar} \big)$ \cite{FS}.
\end{example}

\vspace*{1cm}
\noindent
\parbox[t]{20em}
{\scriptsize{
\noindent
Giuseppe Dito\\
Institut de Math{\'e}matiques de Bourgogne\\
Universit{\'e} de Bourgogne\\
B.P. 47870,
21078 Dijon Cedex,
France\\
mailto: giuseppe.dito@u-bourgogne.fr\\
http://www.u-bourgogne.fr/monge/g.dito/}}
\qquad
\parbox[t]{15em}
{\scriptsize{
Pierre Schapira\\
Institut de Math{\'e}matiques\\
Universit{\'e} Pierre et Marie Curie\\
175, rue du Chevaleret,
75013 Paris, France\\
mailto: schapira@math.jussieu.fr\\
http://www.math.jussieu.fr/$\sim$schapira/}}

\end{document}